\documentclass[runningheads]{llncs}
\usepackage{graphicx}
\usepackage{amsmath}
\usepackage{url}
\usepackage{amssymb}
\usepackage{listings}
\usepackage{lstlangcoq}
\newcommand\Coquelicot{\lstinline{Coquelicot}}
\newcommand\mathcomp{\lstinline{MathComp}}

\begin{document}
\title{Formalization of Asymptotic Convergence for Stationary Iterative Methods}
%
%
\author{Mohit Tekriwal\thanks{currently affiliated with the Lawrence Livermore National Laboratory, USA.}\and
Joshua Miller \and
Jean-Baptiste Jeannin}
\authorrunning{Mohit et al.}
%
\institute{University of Michigan, Ann Arbor MI 48109, USA\\
\email{\{tmohit, joshmi, jeannin\}@umich.edu}}
\maketitle              
\begin{abstract}
Solutions to differential equations, which are used to model physical systems, are computed numerically by solving a set of discretized equations. This set of discretized equations is reduced to a large linear system, whose solution is typically found using an iterative solver. We start with an initial guess, $x_0$, and iterate the algorithm to obtain a sequence of solution vectors, $x_k$, which are approximations to the exact solution of the linear system, $x$. The iterative algorithm is said to converge to $x$, in the field of reals, if and only if $x_k$ converges to $x$ in the limit of $k \to \infty$.

In this paper, we formally prove the asymptotic convergence of a particular class of iterative methods called the \emph{stationary iterative methods}, in the Coq theorem prover. We formalize the necessary and sufficient conditions required for the \emph{iterative convergence}, and extend this result to two classical iterative methods: the Gauss--Seidel method and the Jacobi method. For the Gauss--Seidel method, we also formalize a set of \emph{easily testable conditions} for iterative convergence, called the \emph{Reich theorem}, for a particular matrix structure, and apply this on a model problem of the one-dimensional heat equation. We also apply the main theorem of iterative convergence to prove convergence of the Jacobi method on the model problem.  


\end{abstract}
\section{Introduction}
Solutions to differential equations are often obtained numerically, which often involves solving a large linear system, $Ax = b$. This system is obtained after discretizing a differential equation in a finite computational domain to obtain a set of discretized equations. Direct methods to solve this linear system, such as Gaussian elimination, usually involve matrix inversion, which is computationally expensive with computational complexity in $\mathcal{O}(N^3)$, where $N$ is the dimension of the linear system. Therefore, low-cost methods such as \emph{iterative methods}~\cite{saad2003iterative}, which have an average complexity in $\mathcal{O}(N^2)$, are used to obtain an \emph{approximate} solution of the linear system. 

The goal of an iterative method is to build a sequence of approximations of the \emph{true numerical solution}, which is defined as $x \overset{\Delta}{=}A^{-1}b$. One starts with an initial guess vector $x_0$, and builds a sequence of approximate solutions: $\{x_1, x_2, \hdots, x_{k-1}, x_k \}$ for $k$ iterations, with the hope that $x_k$ is close to $x$. The distance between $x_k$ and $x$ is called the \emph{iterative convergence error}. To ensure the asymptotic convergence of these approximate solutions to the true numerical solution, we need to bound the iterative convergence error, and further show that this error decreases as we increase the number of iterations.

Many general purpose ordinary differential equation (ODE) solvers use some kind of iterative method to solve linear systems. For instance, ODEPACK~\cite{hindmarsh1983odepack}, which is a collection of FORTRAN solvers for initial value problem for ODEs, uses iterative (preconditioned Krylov) methods instead of direct methods for solving linear systems. 
Another widely used suite of ODE solvers is the SUNDIALS~\cite{hindmarsh2005sundials}. 
SUNDIALS has support for a variety of direct and Krylov iterative methods for solving the system of linear equations.
SUNDIALS solvers are used by the mixed finite element (MFEM) package for solving nonlinear algebraic systems and by NASA for spacecraft trajectory simulation~\cite{hindmarsh2005sundials}.
Because those iterative methods are widely used, it is important to obtain formal guarantees for the convergence of iterative solutions to the ``true'' solution of differential equations. 
In this work we use the Coq theorem prover to formalize the convergence guarantees for a class of iterative methods called the \emph{Stationary iterative methods}. The choice of stationary iterative methods for formalization is motivated by the fact that this class of methods is used as building blocks for more complicated iterative solvers like Krylov subspace and conjugate gradient methods~\cite{saad2003iterative}.

\vspace{-0.8em}
\paragraph{Contributions:}
We provide an overview of the \emph{stationary iterative methods} in Section~\ref{sec:overview}, followed by formalization of a generalized iterative convergence theorem in Coq, and its specialization to two classical iterative methods. Overall, this work 
\footnote{Our Coq formalization is available at \url{https://github.com/mohittkr/iterative_convergence.git}}
makes the following contributions:
\begin{itemize}
     \item We provide a formalization of the necessary and sufficient conditions for iterative convergence in Coq in Section~\ref{sec:iter_theorem_proof};
    \item We formalize a set of easily testable conditions for convergence of the Gauss--Seidel classical iterative method for a specific matrix structure in Section~\ref{sec:gauss_reich} and prove convergence on a model problem in Section~\ref{sec:gauss_seidel_model}; 
    \item We then apply the generalized iterative convergence theorem to an example of the Jacobi iteration, another classical iterative method, to prove its convergence in Section~\ref{sec:jacobi_proof};
    \item We develop libraries for dealing with complex matrices and vectors, and formalize $\ell^2$ norm of a matrix and its spectral properties.
\end{itemize}
All of the above formalization has been done in the field of real numbers. 

\section{Overview of Stationary Iterative Methods}\label{sec:overview}
In this section, we provide an overview of the stationary iterative methods adapted from the textbook~\cite{saad2003iterative}.

\begin{figure}[ht]
    \centering
    \includegraphics[scale =0.22]{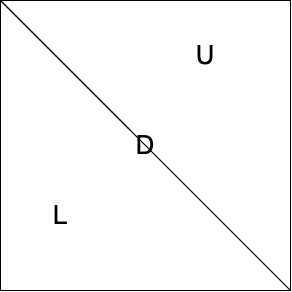}
    \caption{Initial partitioning of matrix $A = L + D + U$. $L$ is the strictly lower triangular matrix. $D$ is the diagonal matrix. $U$ is the strictly upper triangular matrix.}
    \label{fig:iterative_splitting}
\end{figure}

Let $x$ be the true numerical solution, or the direct solution obtained by inverting the linear system $Ax = b$ as
\begin{equation}\label{direct_solution}
    x  \overset{\Delta}{=} A^{-1}b
\end{equation}
Here, the \emph{coefficient} matrix $A \in \mathbb{R}^{n \times n}$ and the \emph{right hand side} vector $b \in \mathbb{R}^n$ are known to us and we are computing the unknown vector $x \in \mathbb{R}^n$. We assume that the matrix $A$ is non-singular. Thus, there exists a unique solution $x$ of the linear system $Ax = b$. For any iterative algorithm, we start with an initial guess vector $x_0$ and obtain a sequence of numerical solutions which are an approximation of the solution $x$. Let $x_k$ be the iterative solution obtained after $k$ iterations obtained by solving the iterative system
\begin{equation}\label{eqn:iterative_system}
    M x_k + N x_{k-1} = b
\end{equation}
for some choice of initial solution vector $x_0$. The vector $x_{k-1}$ is the iterative solution obtained after $k-1$ iterations. At the $k^{th}$ iteration step, $x_{k-1}$ is known to us. The matrices $M$ and $N$ are obtained by splitting (\emph{regular splitting}~\cite{saad2003iterative}) the original coefficient matrix $A$ such that $M$ is easily invertible. Therefore,
\begin{equation}\label{eqn:matrix_split}
    A = M + N.
\end{equation}
The choice of matrices $M$ and $N$ define the choice of an iterative method. For instance, if we choose $M$ to be the lower triangular entries of $A$ and $N$ to be the strictly upper triangular entries of $A$, we get the Gauss--Seidel iterative method~\cite{saad2003iterative}. Thus, for the Gauss--Seidel method, $M = L+D$, and $N = U$, where $L$, $D$, and $U$ are illustrated in Figure~\ref{fig:iterative_splitting}. 
If we choose $M$ to be the diagonal entries of matrix $A$ and $N$ to be the strictly lower and upper triangular entries of $A$, we obtain the Jacobi method~\cite{saad2003iterative}. 
Thus, for the Jacobi method, $M = D$, and $N = L + U$.
Therefore, the matrices $M$ and $N$ are also known to us based on the choice of an iterative method. The right hand vector $b$ is also known to us. Thus, the unknown solution vector $x_k$ at $k^{th}$ step is obtained by rearranging terms in the iterative system~(\ref{eqn:iterative_system}) as
\begin{equation}\label{eqn:iterative_solution}
    x_k = (-M^{-1}N)x_{k-1} + M^{-1}b  
\end{equation}
The quantity $(-M^{-1}N)$ in equation~(\ref{eqn:iterative_solution}) is called an \emph{iterative matrix} and we will denote it as $S$. Therefore,
\begin{equation}\label{eqn:iterative_matrix}
   S \overset{\Delta}{=} -M^{-1}N
\end{equation}
The iterative convergence error after $k$ iterations is defined as
\begin{equation}\label{eqn:error_recurrence}
    e_k^{iterative} \overset{\Delta}{=}
    x_k - x = S^ke_0; \quad e_o \overset{\Delta}{=} x_o - x
\end{equation}
The iterative solution $x_k$ is said to converge to $x$ in the field of reals if and only if
\begin{equation}
    \lim_{k \to \infty} ||e_k^{iterative}|| = \lim_{k \to \infty} ||x_k - x|| = 0
\end{equation}
where $||.||$ denotes a vector norm. In this paper, we will be using the $\ell^2$ vector norm defined as $||x||_2 = \sqrt{\sum_{j=1}^n |x_{\{i\}}|^2}.$

In this work, we consider two concrete instances of stationary iterative methods: the Gauss--Seidel method~\cite{saad2003iterative}, for which the iterative matrix $S_G \overset{\Delta}{=} -M^{-1}N = - (L+D)^{-1}U$, and the Jacobi method~\cite{saad2003iterative}, for which the iterative matrix $S_J \overset{\Delta}{=} -M^{-1}N =  -D^{-1}(A-D) = I-D^{-1}A$.

\vspace{-1em}
\section{Generic iterative convergence theorem in the field of reals}\label{sec:iter_theorem}
\vspace{-0.5em}
The following theorem provides necessary and sufficient conditions for iterative convergence in the field of reals.
\begin{theorem}\label{nec_suf}~\cite{saad2003iterative}
Let an iterative matrix be defined as~(\ref{eqn:iterative_matrix}) for the iterative system~(\ref{eqn:iterative_system}). The sequence of iterative solutions $\{x_k\}$ converges to the direct solution $x$ for all initial values $x_0$, if and only if the spectral radius of the iterative matrix $S = -M^{-1}N$ is less than $1$, i.e.,
\begin{equation*}
    (\forall~x_o,\lim_{k \to \infty}||x_k- x|| = 0) \iff
    \rho(-M^{-1} N)< 1.
\end{equation*}
\end{theorem}
\noindent The spectral radius $\rho$ of a matrix is defined as the maximum eigenvalue in magnitude. We next discuss the proof of Theorem~\ref{nec_suf} followed by its formalization in the Coq proof assistant. It is noteworthy that while such proofs have been discussed in numerical analysis literature, we found several missing pieces during the formalization. Most facts about intermediate steps in the proof were just stated in the numerical analysis literature without a rigorous proof. 
In that regard, a contribution of this work is to provide a clean machine-checked proof of the main theorem and any intermediate lemma or fact that was required to close the proof of Theorem~\ref{nec_suf}.
Because of a lack of space, we provide most of the informal proofs in the Appendix~\ref{thm:lem1_proof} and \ref{thm:lem2_proof}.

To prove Theorem~\ref{nec_suf}, we first need to obtain a recurrence relation for the iterative convergence error at $k^{th}$ step in terms of the initial iteration error $(x_0 - x)$. Therefore,
\begin{proof}[Proof of Theorem~\ref{nec_suf}]
\begin{align}
    x_{k}-x &= -M^{-1}N x_{k-1}+M^{-1}b -x \nonumber \\
    &= -M^{-1}Nx_{k-1}+ M^{-1}(Ax) -x \quad [\text{Since,}~ Ax\overset{\Delta}{=} b] \nonumber \\
    &= -M^{-1}Nx_{k-1}+M^{-1}(M+N)x -x~[\text{Since, }M+N = A] \nonumber \\
    &=-M^{-1}Nx_{k-1}+M^{-1}Mx +M^{-1}Nx -x \nonumber \\
    &=-M^{-1}N(x_{k-1}-x) \quad [\text{Since, }M^{-1}M \overset{\Delta}{=} I]\nonumber
\end{align}
Taking the norm of the vector on both sides, the iterative convergence error at the $k^{th}$ step can be written in terms of the iterative convergence error at $(k-1)$ step as
\begin{equation}\label{eqn:iter1}
    ||x_k - x|| = || (-M^{-1}N) (x_{k-1} -x) ||
\end{equation}
Since, the system is linear, equation~(\ref{eqn:iter1}) can be written in terms of the initial iteration error as
\begin{equation}\label{eqn:iter2}
    ||x_k - x|| = || (-M^{-1}N)^k (x_0 - x)||
\end{equation}
Taking limits of the vector norms on both sides of equation~(\ref{eqn:iter2}),
\begin{equation}\label{eqn:iter3}
    \lim_{k \to \infty} ||x_k - x|| = \lim_{k \to \infty} ||(-M^{-1}N)^k (x_0 - x)||
\end{equation}
If $x_0 = x$, the iterative convergence error is zero trivially. The case $x_0 \neq x$ is interesting and we can prove Theorem~\ref{nec_suf} by splitting it into two lemmas
\begin{lemma}\label{lem:subproof_1}
For given matrices $M \in \mathbb{R}^{n \times n}$ and
$N \in \mathbb{R}^{n \times n}$ respecting the regular splitting, i.e., $A = M + N$,
the sequence of iterative solutions $\{x_k\}$ converges to $x$ for any given initial vector $x_0$ if and only if the $\ell^2$ matrix norm of the iterated product of the iteration matrix, $(-M^{-1}N)^k$ approaches zero as $k \to \infty$, i.e.,
\begin{equation}\label{eqn:subproof_1}
(\forall x_0,  \lim_{k \to \infty} ||(-M^{-1}N)^k (x_0 - x)|| = 0) \iff 
\lim_{k \to \infty} ||(-M^{-1}N)^k|| = 0
\end{equation}
\end{lemma}

\begin{lemma}\label{lem:subproof_2}
For given matrices $M \in \mathbb{R}^{n \times n}$ and
$N \in \mathbb{R}^{n \times n}$ respecting the regular splitting, i.e., $A = M + N$,
the $\ell^2$ matrix norm of the iterated product of the iteration matrix, $(-M^{-1}N)^k$ approaches zero as $k \to \infty$ if and only if the spectral radius of the iteration matrix is less than $1$, i.e., 
\begin{equation}\label{eqn:subproof_2}
       \lim_{k \to \infty} ||(-M^{-1}N)^k|| = 0 \iff \rho(-M^{-1}N) < 1 
\end{equation}
\end{lemma}
\noindent Therefore, by composing proofs of the Lemma~\ref{lem:subproof_1} (see Appendix~\ref{thm:lem1_proof}) 
and the Lemma~\ref{lem:subproof_2} (see Appendix~\ref{thm:lem2_proof}) 
we close the proof of the Theorem~\ref{nec_suf}.
\end{proof}



\section{Formalization of the the generalized iterative convergence theorem (Theorem~\ref{nec_suf}) in Coq}\label{sec:iter_theorem_proof}

One of the main challenges that we encountered during our formalization was the lack of theories for matrix and vector norm and generic properties about complex vectors and matrices. We will first discuss the formalization of these properties in Coq, followed by their adoption in our formalization on iterative convergence.

\subsection{Formalizing properties of complex matrices and vectors}

The \texttt{complex} theory in the \texttt{real\_closed}~\cite{cohen:hal-00671809} library in \mathcomp~ defines complex numbers and basic operations on them. They define complex numbers as a real closed field, thereby allowing us to instantiate a generic field with a complex field. This was useful when we used the \texttt{eigenvalue} definition from \mathcomp~ matrix algebra library and the canonical forms library by Cano~et~al~\cite{cano:hal-00779376}. However, since the
basic properties like modulus of a complex number, conjugates, properties of complex matrices and vectors were lacking, we added them in our formalization\footnote{The theory about complex modulus and norms has been added in most recent developments of \mathcomp~after our discussion with the developers.
We were pointed to the development of matrix norms in the CoqQ project (\url{https://github.com/coq-quantum/CoqQ/blob/main/src/mxnorm.v}), which was done concurrently with our development.
}.
We define the modulus of a complex number as
\begin{lstlisting}
Definition C_mod (x: complex R):= sqrt ( (Re x)^+2 + (Im x)^+2).
\end{lstlisting}
\texttt{Re x} and \texttt{Im x} denote the real and imaginary part of \texttt{x}, respectively. We proved some basic properties of the modulus, which we enumerate in Table~\ref{tab:modulus}.~
\begin{table*}[htbp]
    \centering
    \begin{tabular}{|c|l|}
    \hline
    Mathematical properties &
    Formalization in Coq \\
    \hline
    $|| 0 || = 0$     & 
    \begin{lstlisting}
    Lemma C_mod_0: C_mod 0 = 0%Re.
    \end{lstlisting}
    \\
    \hline
    $0 \leq ||x||$   & 
    \begin{lstlisting}
    Lemma C_mod_ge_0: forall (x: complex R), 
        (0<= C_mod x)%Re.
    \end{lstlisting}\\
    \hline
    $||x y|| = ||x||~ ||y||$     & 
    \begin{lstlisting}
    Lemma C_mod_prod: forall (x y: complex R),
        C_mod (x * y) = C_mod x * C_mod y.
    \end{lstlisting} \\
    \hline 
    $||\frac{x}{y}|| =  \frac{||x||}{||y||}, y \neq 0$ &
    \begin{lstlisting}
    Lemma C_mod_div: forall (x y: complex R), y <> 0 -> 
        C_mod (x / y) = (C_mod x) / (C_mod y).
    \end{lstlisting}\\
    \hline
    $||x|| \neq 0,\quad \text{if. } x \neq 0$ &
    \begin{lstlisting}
    Lemma C_mod_not_zero: forall (x: complex R), 
        x <> 0 -> C_mod x <> 0.
    \end{lstlisting}\\
    \hline
    $||1|| = 1$ &
    \begin{lstlisting}
    Lemma C_mod_1: C_mod 1 = 1.
    \end{lstlisting}\\
    \hline
    $||x^n|| = ||x||^n$ &
    \begin{lstlisting}
    Lemma C_mod_pow: forall (x: complex R) (n:nat),  
        C_mod (x^+ n) = (C_mod x)^+n.
    \end{lstlisting}\\
    \hline
    $||x + y|| \leq ||x|| + ||y||$ &
    \begin{lstlisting}
    Lemma C_mod_add_leq : forall (a b: complex R), 
    C_mod (a + b) <=  C_mod a + C_mod b.
    \end{lstlisting}\\
    \hline
    $||\frac{1}{x}|| = \frac{1}{||x||}, \quad \text{if  } x \neq 0$ &
    \begin{lstlisting}
    Lemma C_mod_inv : forall x : complex R, x <> 0 -> 
        C_mod (invc x) = Rinv (C_mod x).
    \end{lstlisting}\\
    \hline
    $||xy||^2 = ||x||^2 ||y||^2$ &
    \begin{lstlisting}
    Lemma C_mod_sqr: forall (x y : complex R), 
        Rsqr (C_mod (x * y)) = 
        (Rsqr (C_mod x)) * (Rsqr (C_mod y)).
    \end{lstlisting}\\
    \hline
    $||-x|| = ||x||$ &
    \begin{lstlisting}
    Lemma C_mod_minus_x: forall (x: complex R), 
        C_mod (-x) = C_mod x.
    \end{lstlisting}\\
    \hline
    $||\sum_{j=0}^n u(j)|| \leq 
    \sum_{j=0}^n ||u(j)||
    $ &
    \begin{lstlisting}
    Lemma C_mod_sum_rel: 
    forall (n:nat) (u : 'I_n.+1 -> (complex R)), 
    (C_mod ($\sum_j$ (u j))) <= $\sum_j$ ((C_mod (u j))).
    \end{lstlisting}\\
    \hline
    \end{tabular}
    \caption{Formalization of properties of complex modulus in Coq}
    \label{tab:modulus}
    \vspace{-8pt}
\end{table*}
We also found that some theory on complex conjugates were missing in the \mathcomp~library. The Table~\ref{tab:conjugate}~
lists the missing formalization of complex conjugates that we added for this formalization. We also formalize the properties of a conjugate matrix like idempotent property, scaling of a complex matrix, conjugate transpose of matrix multiplication, etc., which we do not list here for brevity.
\begin{table*}[htbp]
    \centering
    \begin{tabular}{|c|l|}
    \hline
    Mathematical properties & Formalization in Coq\\
    \hline
    $\overline{xy} = \bar{x} \; \bar{y}$    & 
    \begin{lstlisting}
    Lemma Cconj_prod: forall (x y: complex R),
        conjc (x*y)%C = (conjc x * conjc y)%C.
    \end{lstlisting} \\
    \hline
    $\overline{x + y} = \bar{x} + \bar{y} $   & 
    \begin{lstlisting}
    Lemma Cconj_add: forall (x y: complex R), 
        conjc (x+y) = conjc x + conjc y.
    \end{lstlisting}\\
    \hline
    $||x|| = ||\bar{x}||$ &
    \begin{lstlisting}
    Lemma Cconjc_mod: forall (a: complex R),
        C_mod a = C_mod (conjc a).
    \end{lstlisting}\\
    \hline
    $x = \bar{\bar{x}}$ &
    \begin{lstlisting}
    Lemma conj_of_conj_C: forall (x: complex R), x = conjc (conjc x).
    \end{lstlisting}\\
    \hline
    $\bar{x} x = ||x||^2$ &
    \begin{lstlisting}
    Lemma conj_prod: forall (x:complex R),
        ((conjc x)*x)%C = RtoC (Rsqr (C_mod x)).
    \end{lstlisting}\textsuperscript{1}\\
    \hline
    $Re[x] + Re[\bar{x}] = 2 Re[x]$ &
    \begin{lstlisting}
    Lemma Re_conjc_add: forall (x: complex R), 
        Re x + Re (conjc x) = 2 * (Re x).
    \end{lstlisting}\\
    \hline
    $\overline{\sum_{j=0}^n f(i)} = \sum_{j=0}^n \overline{f(i)}$ &
    \begin{lstlisting}
    Lemma Cconj_sum: forall (p:nat) (x: 'I_p -> complex R),
    conjc ($\sum_{j < p}$  x j)= $\sum_{j < p}$ conjc (x j).
    \end{lstlisting}\\
    \hline
    \end{tabular}
    \caption{Formalization of properties of complex conjugates in Coq. \textsuperscript{1}Here, \lstinline{RtoC} is a coercion from reals to complex. }
    \label{tab:conjugate}
    \vspace{-14pt}
\end{table*}


\subsection{Formalization of vector and matrix norms}
Another missing piece in the existing linear algebra theory in \mathcomp~was the norm of a vector and a matrix.
In this work, we formalize the $\ell^2$-norm of a vector and its induced matrix norm. In Coq, we define the $\ell^2$-norm of a matrix as
\begin{lstlisting}
Definition matrix_norm (n:nat) (A: 'M[complex R]_n.+1) := Lub_Rbar (fun x =>
    exists v: 'cV[complex R]_n.+1, v != 0 /\ x = (vec_norm_C  (A *m v))/ (vec_norm_C v)).
\end{lstlisting}
where \texttt{vec\_norm\_C} is the $\ell^2$-norm of a complex vector, which we define in Coq as
\begin{lstlisting}
Definition vec_norm_C (n:nat) (x: 'cV[complex R]_n.+1):= sqrt ($\sum_l$ (C_mod x l 0)^2)
\end{lstlisting}

\noindent The definition \texttt{Lub\_Rbar} is the least upper bound and is already defined in the \texttt{Coquelicot}~\cite{boldo2015coquelicot} library. Mathematically, \texttt{matrix\_norm} formalizes the following definition of a matrix norm
$||A||_i = \sup_{x \neq 0}\frac{||Ax||}{||x||}$,
for a given vector norm $||.||$, which in this case is the $\ell^2$ vector norm. 
In Table~\ref{tab:vector_norm_table} and Table~\ref{tab:matrix_norm_table} 
we enumerate the properties of vector and matrix norms that we formalized.

\begin{table*}[htbp]
    \centering
    \begin{tabular}{|c|l|}
    \hline
    Mathematical properties & Formalization in Coq\\
    \hline
    $0 \leq ||v||$    & \begin{lstlisting}
    Lemma vec_norm_C_ge_0: forall (n:nat)
    (v: 'cV[complex R]_n.+1),  0<= vec_norm_C v.
    \end{lstlisting} \\
    \hline
    $||av|| = |a|\; ||v||, \quad \text{$a$ is scalar} $ & 
    \begin{lstlisting}
    Lemma ei_vec_ei_compat:
    forall (n:nat) (x:complex R) (v: 'cV[complex R]_n.+1), 
     vec_norm_C (scal_vec_C x v) =
        C_mod x * vec_norm_C v.
    \end{lstlisting}\\
    \hline
    $||v_1 + v_2|| \leq ||v_1|| + ||v_2||$ &
    \begin{lstlisting}
    Lemma vec_norm_add_le: 
    forall (n:nat) (v1 v2 : 'cV[complex R]_n.+1),
    vec_norm_C (v1 + v2) <= 
        vec_norm_C v1 + vec_norm_C v2.
    \end{lstlisting}\\
    \hline
    $v \neq 0 \implies ||v|| \neq 0$\textsuperscript{1} &
    \begin{lstlisting}
    Lemma non_zero_vec_norm: 
    forall (n:nat) (v: 'cV[complex R]_n.+1),
     vec_not_zero v -> vec_norm_C v <> 0.
    \end{lstlisting}\\
    \hline
    \end{tabular}
    \caption{Formalization of properties of vector norm in Coq. \textsuperscript{1} \lstinline{vec_not_zero} is Coq's definition of $v \neq 0$.}
    \label{tab:vector_norm_table}
\vspace{-20pt}
\end{table*}

\begin{table*}[htbp]
    \centering
    \begin{tabular}{|c|l|}
    \hline
    Mathematical properties & Formalization in Coq \\
    \hline
    $0 \leq ||A||_i$    & \begin{lstlisting}
    Lemma matrix_norm_ge_0:
    forall (n:nat) (A: 'M[complex R]_n.+1), 
        0 <= matrix_norm A.
    \end{lstlisting} \\
    \hline
    $||Ax|| \leq ||A||_i ||x||, \quad x \neq 0 $\textsuperscript{1}  & 
    \begin{lstlisting}
    Lemma matrix_norm_compat: 
    forall (n:nat) (x: 'cV[complex R]_n.+1) 
    (A: 'M[complex R]_n.+1), x != 0 -> 
    vec_norm_C (mulmx A x) <= 
        (matrix_norm A) * vec_norm_C x.
    \end{lstlisting}\\
    \hline
    $||AB||_i \leq ||A||_i ||B||_i$ &
    \begin{lstlisting}
    Lemma matrix_norm_prod:
    forall (n:nat) (A B: 'M[complex R]_n.+1),
    matrix_norm (A *m B) <= 
        (matrix_norm A) * (matrix_norm B).
    \end{lstlisting}\\
    \hline
    $0 \leq ||A||_i \leq ||A||_F$\textsuperscript{2} &
    \begin{lstlisting}
    Lemma mat_2_norm_F_norm_compat: 
    forall (n:nat) (A: 'M[complex R]_n.+1),
    0 <= matrix_norm A <= mat_norm A.
    \end{lstlisting}\\
    \hline
    \end{tabular}
    \caption{Formalization of properties of matrix norm in Coq. \\ \textsuperscript{1}Here, $x$ is a vector and the relation proves compatibilty relation between a matrix norm and its induced vector norm. \\
    \textsuperscript{2}Here, $||A||_F$ is the Frobenius norm and we prove that the 2-norm of a matrix is bounded above by the Frobenius matrix norm.
    The Frobenius norm of a matrix is defined as $||A||_F = \sqrt{\sum_{j=1}^n \sum_{j=1}^n |A_{ij}|^2}$}   \label{tab:matrix_norm_table}
    \vspace{-14pt}
\end{table*}
\noindent An important point to note here is that since we are using the \texttt{Coquelicot} definition of an extended real line, \texttt{Rbar}, coercion of a quantity of type \texttt{Rbar} to real requires us to prove finiteness of that quantity. We therefore have to prove that the matrix norm is finite, which we state as the following lemma in Coq
\begin{lstlisting}
Lemma matrix_norm_is_finite: forall (n:nat) (A: 'M[complex R]_n.+1), is_finite (matrix_norm A).
\end{lstlisting}
Since we prove asymptotic convergence of component-wise limit of the elements of a Jordan matrix, we have to work with the Frobenius norm of a matrix, which we define in Coq as
\begin{lstlisting}
Definition mat_norm (n:nat) (A: 'M[complex R]_n.+1) : R := sqrt ($\sum_i$ $\sum_j$ (C_mod (A i j))^2)
\end{lstlisting}

\noindent We will next discuss
the formalization of the Theorem~\ref{nec_suf}.
\vspace{-1em}
\subsection{Formalization of the Theorem~\ref{nec_suf} in Coq}
We state Theorem~\ref{nec_suf} in Coq as follows:
\begin{lstlisting}
Theorem iter_convergence: forall (n:nat) (A: 'M[R]_n.+1) (b: 'cV[R]_n.+1) (M N : 'M[R]_n.+1), 
A \in unitmx -> M \in unitmx -> A = M + N ->
let x := (A^-1) *m b in
(let S_mat := RtoC_mat (- ( M^-1 *m N)) in 
(forall (i: 'I_n.+1), (C_mod (lambda S_mat i) < 1))) <->
(forall x0: 'cV[R]_n.+1, is_lim_seq (fun k:nat => vec_norm ((X_m k.+1 x0 b M N) - x)) 0).
\end{lstlisting}
The theorem \lstinline{iter_convergence} states that if the square matrix $A \in \mathbb{R}^{n+1 \times n+1}$ is invertible, which is formalized by the condition \lstinline{A \in unitmx} and if the sub-matrices $M, N \in \mathbb{R}^{n+1 \times n+1}$ respect the regular splitting property, which is formalized by $A = M+N$, and if $M$ is invertible, which is formalized by \lstinline{M \in unitmx}, then $\lim_{k \to \infty} ||x_{k+1} - x|| = 0$ if and only if $\forall~i, |\lambda_i (S)| < 1$. Here, $|.|$ is the complex modulus of the eigenvalues of the iteration matrix $S$ (\lstinline{S_mat} as defined in the let statement of the theorem \lstinline{iter_convergence}), which is defined as in equation~(\ref{eqn:iterative_matrix}).
The \lstinline{is_lim_seq} predicate from the \Coquelicot~library~\cite{boldo2015coquelicot}  defines limit of a sequence.
The reason for having the dimension of $A$ as $(n+1) \times (n+1)$ is that the natural numbers start from $0$ in Coq and a square matrix of size $0$ does not type check. Thus, instead of having $0 < n$ as pre-condition in the theorem explicitly, we specify this implicitly in the dimension of $A$. This change also makes our development compatible with the Jordan canonical formalization by Cano~et~al~\cite{cano:hal-00779376}, who need this constraint on the dimension to define and work with block matrix.

Since we deal with a generic case where a real matrix is allowed to have complex eigenvalues and eigenvectors, we need to transform the real iteration matrix $S$ to a complex matrix, so as to be consistent with types. Thus, given a real matrix $A$, \lstinline{RtoC_mat} transforms a real entry $A_{ij}: \mathbb{R}$ to a complex number $\tilde{A}_{ij}:\mathbb{C}:= (A_{ij} +i* 0)$. 
An important point to note here is that we do not compute the \emph{true numerical solution} $x$ explicitly, but define it as $x  \overset{\Delta}{=} A^{-1}b$ in the let binding of the theorem statement -- \lstinline{let x := (A^-1) *m b}.
We define the iterative solution after $k$ steps, $x_k$ from the iterative system~(\ref{eqn:iterative_system}) using the \lstinline{Fixpoint} operator in Coq, which lets us define the recurrence relation~(\ref{eqn:iterative_solution})
\begin{lstlisting}
Fixpoint X_m (k n:nat) (x0 b: 'cV[R]_n.+1) (M N: 'M[R]_n.+1) : 'cV[R]_n.+1:=
match k with
| O => x0
| S p => ((- ((M^-1) *m N)) *m (X_m p x0 b M N)) + ((M^-1) *m b)
end.
\end{lstlisting}

\vspace{-1.8em}
\subsubsection{Formalization of eigenvalues}
One of the main issues we faced was coming up with a scalar definition of eigenvalues, which satisfies the characteristic definition $Av = \lambda v$. The \mathcomp~library defines a non-computable definition of eigenvalues, \lstinline{eigenvalue}, which is a predicate stating that $\lambda$ is an eigenvalue of a matrix $A$, if the eigenspace corresponding to $\lambda$ is non-zero. For the scalar definition of eigenvalue, we took an inspiration from the Jordan canonical forms formalization by Guillaume Cano and Maxime  D{\'e}n{\`e}s~\cite{cano:hal-00779376}. They define a Jordan form, whose diagonal elements are the characteristic polynomials of the Smith Normal form of a matrix $A$. Since the diagonal elements of a Jordan form are the eigenvalues of that matrix, we then define a sequence of eigenvalues from these diagonal entries of the Jordan matrix as
\begin{lstlisting}
Definition lambda_seq (n: nat) (A: 'M[complex R]_n.+1) :=
let sizes:= size_sum [seq x.2.-1 | x <- root_seq_poly (invariant_factors A)] in 
[seq (Jordan_form A) i i | i <- enum 'I_sizes.+1].
\end{lstlisting}
where \lstinline{root_seq_poly p} returns a sequence of pair of roots and its multiplicity, of the polynomial \lstinline{p}. The \lstinline{invariant_factors} are the polynomials in the diagonal of the Smith Normal form of a matrix. In this case, the sequence contains the pair of eigenvalues of matrix $A$ and its multiplicity. The $i^{th}$ eigenvalue of $A$ is then defined as the $i^{th}$ component of the sequence of eigenvalues \lstinline{lambda_seq}.
\begin{lstlisting}
Definition lambda (n: nat) (A: 'M[complex R]_n.+1) (i: 'I_n.+1) := 
    (@nth _  0%C (lambda_seq A) i).
\end{lstlisting}
To take full advantage of the lemmas describing eigenvalues and eigenvectors as defined in \mathcomp, we had to relate the definition of eigenvalue, \texttt{lambda}, and the one defined in \mathcomp~. 
The lemma \lstinline{Jordan_ii_is_eigen} asserts that \lstinline{lambda} satisfies the predicate \lstinline{eigenvalue}, and is indeed an eigenvalue of a matrix \lstinline{A}.
\begin{lstlisting}
Lemma Jordan_ii_is_eigen: forall (n: nat) (A: 'M[complex R]_n.+1),
  forall  (i: 'I_n.+1), @eigenvalue (complex_fieldType _) n.+1 A (@nth _  0%C (lambda_seq A) i).
\end{lstlisting}
Here, \texttt{size\_sum} is the sum of the algebraic multiplicities of the eigenvalues and equals the total size of the matrix $n$. We prove this fact using the following lemma statement in Coq
\begin{lstlisting}[language = Coq]
Lemma total_eigen_val: forall (n:nat) (A: 'M[complex R]_n.+1),
(size_sum [seq x.2.-1 | x <- root_seq_poly (invariant_factors A)]).+1 = n.+1.
\end{lstlisting}
The lemma \texttt{total\_eigen\_val} helps us get around the dimension constraint imposed by the design of the Jordan form of a matrix $A$.

\vspace{-1em}
\subsubsection{Formalization of the ratio test}
As discussed in the informal proof in Appendix~\ref{thm:lem2_proof}, to prove sufficiency condition for iterative convergence, we had to formalize the ratio test for convergence of sequences which was missing in the existing Coq libraries.
In Coq, we state the ratio test (Lemma~\ref{ratio_test} in Appendix~\ref{thm:lem2_proof}) as: 
\begin{lstlisting}
Lemma ratio_test: forall (a: nat -> R) (L:R), (0 < L /\ L < 1) -> 
(forall n:nat, (0 < a n)) -> 
(is_lim_seq ( fun n:nat =>  ((a (n.+1))/(a n))) L) ->
is_lim_seq (fun n: nat => a n) 0.
\end{lstlisting}
\vspace{-1em}
\subsubsection{Proving convergence of the Jordan block matrix}
To prove that each element of the Jordan block matrix $J^k$ (see Appendix~\ref{thm:lem2_proof} 
for more details) converges to zero, i.e.,
\begin{equation*}
     \forall i \; j, \lim_{k \to \infty} |(J^k)_{(i,j)}|^2 = 0
\end{equation*}
where 
\begin{equation*}
    J^k = 
    \begin{bmatrix}
    J_{m_1}^k(\lambda_1) & 0 & 0 & \hdots & 0\\
    0 & J_{m_2}^k(\lambda_2) & 0 & \hdots & 0 \\
    \vdots & \hdots & \ddots & \hdots & \vdots \\
    0 & \hdots & 0 & J_{m_{s-1}}^k(\lambda_{s-1}) & 0\\
    0 & \hdots & \hdots & 0 & J_{m_s}^k(\lambda_{s})
    \end{bmatrix}
\end{equation*}
and
\begin{equation*}
    J_{m_i}^k(\lambda_i) =
    \begin{bmatrix}
    \lambda_i^k & {k \choose 1}\lambda_i^{k-1} & {k \choose 2}\lambda_i^{k-2} & \hdots & {k \choose {m_i-1}}\lambda_i ^{k - m_i +1} \\
    0 &  \lambda_i^k & {k \choose 1}\lambda_i^{k-1} & \hdots &
    {k \choose {m_i-2}}\lambda_i ^{k - m_i +2}\\
    \vdots & \vdots & \ddots & \ddots & \vdots\\
    0 & 0 & \hdots & \lambda_i^k & {k \choose 1}\lambda_i^{k-1} \\
    0 & 0 & \hdots & 0 & \lambda_i^k
    \end{bmatrix},
\end{equation*} we prove the following lemma.
\begin{lstlisting}
Lemma each_entry_zero_lim: forall (n:nat) (A: 'M[complex R]_n.+1), (i j: 'I_(size_sum sizes).+1),
let sp := root_seq_poly (invariant_factors A) in
let sizes := [seq x0.2.-1 | x0 <- sp] in
(forall i: 'I_(size_sum sizes).+1 , (C_mod (nth 0%C (lambda_seq A) i) < 1) ) ->
is_lim_seq  (fun m: nat => 
let block := (fun n0 i1 : nat => let lambda :=  (nth (0, 0%N) sp i1).1 in 
            \matrix_(i2, j0) (${m.+1}\choose{j0-i2}$ * (lambda ^  (m.+1 - (j0 - i2))) *+ (i2 <= j0)) in 
        (C_mod ((diag_block_mx sizes block) i j))^2) 0.
\end{lstlisting}
The lemma \lstinline{each_entry_zero_lim} states that if the magnitude of each eigenvalue of a matrix $A$ is less than $1$, i.e., $|\lambda_i(A)| < 1, \; \forall i, \; 0 \leq i < N$, then the limit of each term in the expanded Jordan matrix is zero as $k \to \infty$. Here, the block diagonal matrix
\lstinline{diag_block_mx} takes an expanded Jordan block $J_{m_i}^k(\lambda_i), \forall i, 0\leq i < N$ and constructs the Jordan matrix $J^k$. 
We then take a modulus of each entry of $J^k$ and prove that its limit is zero as $k \to \infty$.
A key challenge we faced when proving the lemma \lstinline{each_entry_zero_lim} was extracting each Jordan block of the diagonal block matrix. The diagonal block matrix is defined recursively over a function which takes a block matrix of size $\mu_i$ denoting the algebraic multiplicity of each eigenvalues $\lambda_i$. We had to carefully destruct this definition of diagonal block matrix and extract the Jordan block and the zeros on the off-diagonal entries, which we formalize using the lemma \lstinline{diag_destruct} in Coq. We can then prove the limit on this Jordan block by exploiting its upper triangular structure. 


Limits of the off-diagonal elements can then be trivially proven to be zero since each of those elements are zero. This completes the proof of sufficiency condition for convergence of iterative convergence error.

\vspace{-1em}
\section{Proof of convergence for the Gauss--Seidel method}\label{sec:gauss_reich}
\vspace{-0.5em}
To prove the convergence of a Gauss--Seidel method, we need to prove that the spectral radius of $S_G$ is less than $1$. But computing the eigenvalues of $S_G$ explicitly is almost impossible for a generic matrix. Therefore, we need an easier check to assert that the spectral radius of $S_G$ is indeed less than $1$. The \emph{Reich theorem}~\cite{reich1949convergence} provides a sufficient condition for the spectral radius of $S_G$ to be less than $1$, for a real and symmetric coefficient matrix $A$, with all of the elements in its main diagonal positive. This condition provides a much easier check, which is linear in time, and when layered with the Theorem~\ref{nec_suf}, provides a sufficient condition for convergence of the Gauss--Seidel iteration for a real and symmetric coefficient matrix $A$, with all of the elements in its main diagonal positive.

We next discuss the formalization of the \emph{Reich theorem}, followed by the main convergence theorem for the Gauss--Seidel iteration.
\vspace{-1em}

\subsection{Formalization of easily checkable conditions}
\begin{theorem}[Reich theorem]\label{Reich}~\cite{reich1949convergence}
If A is real, symmetric nth-order matrix with all terms on its main diagonal positive, then a sufficient condition for all the n characteristic roots of $(-M^{-1}N)$ to be smaller than unity in magnitude is that A is positive definite.
\end{theorem}
\noindent From an application point of view, only the sufficiency condition is important. This is because to apply Theorem~\ref{nec_suf}, we only need to know that the magnitude of the eigenvalues are less than 1.  Thus, to prove the convergence of Gauss--Seidel iteration, we first apply Theorem~\ref{nec_suf} to get the eigenvalue condition in the goal and then apply Theorem~\ref{Reich} to complete the proof. Since computing eigenvalues are not very trivial in most cases, the positive definite property of the matrix $A$ provides an easy test for $|\lambda| < 1$ for Gauss--Seidel iteration matrix. 

We next present an informal proof of the Reich Theorem followed by its formalization in Coq
\begin{proof}
Let $z_i$ be the $i^{th}$ characteristic vector of $-(A_{1}^{-1}A_2)$ corresponding to the characteristic root $\mu_{i}$. Then 
\begin{equation}
    -(A_{1}^{-1}A_{2}) z_i = \mu_{i}z_{i}
\end{equation}
Multiplying by $-(\bar{z_{i}}'A_{1})$ on both sides,
\begin{equation}\label{reich_2}
    (-\bar{z_{i}}'A_{1})(-A_{1}^{-1}A_{2})z_i = -\mu_{i}\bar{z_{i}}'A_{1}z_{i}
\end{equation}
where $\bar{z_{i}}'$ is the conjugate transpose of $z_i$ obtained by taking the conjugate of each element of $z_i$ followed by transpose of the vector. Equation~(\ref{reich_2}) then simplifies to:
\begin{equation}
    \bar{z_i}'A_2 z_i = -\mu_i \bar{z_i}'A_1 z_i; \quad [A_1 A_{1}^{-1}=I]
\end{equation}
Consider the bi-linear form, $\bar{z_i}'Az_i$,
\begin{align}
    \bar{z_i}'Az_{i}= \bar{z_i}'(A_1+A_2)z_i &= \bar{z_{i}}'A_1 z_i + \bar{z_i}'A_2 z_i \nonumber \\
    &= \bar{z_{i}}'A_1 z_i - \mu_i \bar{z_i}'A_1 z_i \nonumber \\
    &= (1-\mu_i)\bar{z_{i}}'A_1 z_i  \label{reich_3}
\end{align}
Taking conjugate transpose of equation~(\ref{reich_3}) on both sides,
\begin{equation}\label{reich_4}
    \bar{z_i}'Az_{i} = (1-\bar{\mu_i})\bar{z_{i}}'A_1' z_i
\end{equation}
Let $D$ be the diagonal matrix defined as:
\begin{equation}
    D_{ij}= 
    \begin{cases}
    A_{ij} & \text{ if $i=j$}\\
    0 & \text{if $i \neq j$}
    \end{cases}
\end{equation}
It can be shown that 
\begin{equation}\label{reich_5}
    A_{1}'= D+ A_2
\end{equation}
Substituting equation~(\ref{reich_5}) in equation~(\ref{reich_4}), 
\begin{align}
    \bar{z_i}'Az_{i} &=(1-\bar{\mu_i})\bar{z_{i}}'(D+A_2) z_i \nonumber \\
    &= (1-\bar{\mu_i})\bar{z_{i}}'Dz_i + 
        (1-\bar{\mu_i})\bar{z_{i}}'A_2z_i \nonumber \\
    & = (1-\bar{\mu_i})\bar{z_{i}}'Dz_i +
        \frac{(1-\bar{\mu_i})}{1-\mu_i}\bar{z_i}'Az_i \label{reich_6}
\end{align}
Simplifying equation~(\ref{reich_6}), 
\begin{equation}\label{reich_7}
    (1-\bar{\mu_i} \mu_i)\bar{z_i}'Az_i = (1-\bar{\mu_i})(1-\mu_i)\bar{z_i}'Dz_i
\end{equation}
But, $\bar{\mu_i}\mu_i= |\mu_i|^2$ and 
$(1-\bar{\mu_i})(1-\mu_i)= |1-\mu_i|^2$
Hence, equation~(\ref{reich_7}) simplifies to,
\begin{equation}\label{reich_8}
    (1-|\mu_i|^2)\bar{z_i}'Az_i = (|1-\mu_i|^2) \bar{z_i}'Dz_i
\end{equation}
Since, the diagonal elements of $A$ is positive, i.e., $A_{ii}>0$, $\bar{z_i}'Dz_i$ is positive definite, i.e., $\bar{z_i}'Dz_i >0$.
Since $\bar{z_i}'Dz_i >0$ and $\bar{z_i}'Az_i >0$, $|\mu_i|<1$.
\end{proof}

\vspace{-0.5em}
\subsubsection*{Formalization in Coq:}
We formalize the Theorem~\ref{Reich} in Coq as follows:
\begin{lstlisting}
Theorem Reich_sufficiency: forall (n:nat) (A: 'M[R]_n.+1),
(forall i:'I_n.+1,  A i i > 0) ->
(forall i j:'I_n.+1,   A i j = A j i) ->  
is_positive_definite A -> 
(let S_G := - ( (RtoC_mat (M_G A)^-1) *m (RtoC_mat (N_G A))) in  
    (forall i: 'I_n.+1, C_mod (lambda S_G i) < 1)).
\end{lstlisting} where positive definiteness of a complex matrix $A$ is defined as:
$\forall x \in \mathbb{C}^{n \times 1}, Re\; [x^{*}Ax]>0$.
$x^{*}$ is the complex conjugate transpose of vector $x$ and $Re\; [x^{*}Ax]$ is the real part of the complex scalar $x^{*}Ax$. 
The hypothesis \lstinline{forall i:'I_n.+1,  A i i > 0}
states that all terms in the main diagonal of A are positive. The hypothesis \lstinline{forall i j:'I_n.+1,   A i j = A j i}
states that the matrix A is symmetric. 

\vspace{-1em}
\subsection{Proof of convergence for the Gauss--Seidel method}\label{sec:gauss_seidel_generic}
We then apply the theorem \lstinline{iter_convergence} with \lstinline{Reich_sufficiency} to prove convergence of the Gauss--Seidel iteration method. We formalize the convergence of the Gauss--Seidel iteration method in Coq as

\begin{lstlisting}
Theorem Gauss_Seidel_converges: forall (n:nat) (A: 'M[R]_n.+1) (b: 'cV[R]_n.+1),
let x := (A^-1) *m b in 
A \in unitmx ->
(forall i : 'I_n.+1, 0 < A i i) ->
(forall i j : 'I_n.+1, A i j = A j i) ->
is_positive_definite A ->
(forall x0: 'cV[R]_n.+1, is_lim_seq (fun k:nat =>
        vec_norm ((X_m k.+1 x0 b (M_G A) (N_G A)) - x)) 0).
\end{lstlisting}
Thus, to prove convergence of the Gauss--Seidel method for a real, symmetric matrix with all of its diagonal elements positive, we just need positive-definiteness of $A$, which can be proved by showing positivity of
determinants of all upper-left sub-matrices (principal minors)~\cite{prussing1986principal}. We illustrate this approach using the following example.
\begin{example}\label{example_pd}
To show that the following matrix is positive definite,
\begin{equation*}
 A_{example} =\begin{bmatrix}
    2 & -1 & 0\\
    -1 & 2 & -1\\
    0 & -1 & 2
\end{bmatrix}   
\end{equation*}
we compute the determinants of all possible $k \times k$ upper sub-matrices, i.e.,
\begin{equation*}
    |2| = 2; \quad 
    \begin{vmatrix}
    2 & -1 \\
    -1 & 2
    \end{vmatrix} = 3; \quad 
    \begin{vmatrix}
      2 & -1 & 0\\
    -1 & 2 & -1\\
    0 & -1 & 2  
    \end{vmatrix} = 4.
\end{equation*}
Since the determinant of all $k \times k$ upper sub-matrices are positive, the matrix $A_{example}$ is positive definite. We will be using this approach to show positive definiteness in Section~\ref{sec:gauss_seidel_model}. Note that this approach did not involve computing the eigenvalues of $S_G$ and only relies on the structure of matrix $A_{example}$, thereby making it easy to check for positive definiteness and hence convergence of Gauss--Seidel method, by virtue of the Reich theorem (Theorem~\ref{Reich}). 
\end{example}

\vspace{-1em}
\paragraph{Note on sufficient conditions for convergence of the Jacobi method.} 
The Reich theorem, which we formalized in Coq, provides sufficient conditions for convergence of the Gauss--Seidel method for a symmetric and positive definite matrix. Similar sufficient conditions also exist for convergence of the Jacobi method, which relies on showing \emph{strict row diagonal dominance of the matrix}, or \emph{diagonal dominance and irreducibility of the matrix}~\cite{bagnara1995unified}. A matrix $A$ is said to be strictly diagonally dominant if $|A_{ii}| > \sum_{j \neq i} |A_{ij}|$. This is a much easier check for convergence than computing the eigenvalues of the iterative matrix explicitly. However, we do not formalize this check in this work because the proof of this fact uses the proof of the \emph{ Gers\'gorin-type theorems}~\cite{lancaster1985theory}. The Gers\'gorin-type theorems have not been formalized in Coq, and therefore their adoption in our work would require a separate proof effort for this theorem.

\vspace{-1em}

\section{Model problem}
\vspace{-0.5em}
We apply our convergence theorems on a concrete linear differential equation $\frac{d^2 u}{dx^2} = 1$ for $x \in (0,1)$, with boundary conditions: $u(0) = u(1) = 0$, as a proof of concept. This differential equation is used to model the heat diffusion in a rod, i.e., $1$-D domain. 
We chose a uniform grid with $P$ points in the interior of the $1-$D domain. The grid has a uniform spacing $h$. We will be using a centered difference scheme~\cite{10.1007/978-3-030-76384-8_20} for discretizing the differential equation. Therefore, the difference equation at point $x_i$ in the interior of the $1-$D domain is given by
\begin{equation}\label{diff}
    \frac{-u(x_{i+1})+2u(x_i)-u(x_{i-1})}{h^2}=-1; 
    \quad h=x_{i+1}-x_i=\frac{1}{P+1}
\end{equation}
When we stack the equation~(\ref{diff}) for all points in the interior of the $1-$D domain, we get a linear matrix system
\begin{equation}\label{matrix_A}
    \underbrace{
    \frac{1}{h^2}
    \begin{bmatrix}
    2 & -1 & 0 & 0 & 0 &  \hdots & 0 \\
    -1 & 2 & -1 & 0 & 0 & \hdots & 0\\
    \vdots & & \ddots & \ddots & \ddots  & & \vdots\\
    0 & \hdots & 0 & 0 & -1 & 2 & -1\\
    0 & \hdots & 0 & 0 & 0  & -1 & 2
    \end{bmatrix}}_{A}
    \underbrace{
    \begin{bmatrix}
    u_1 \\
    u_2 \\
    \vdots\\
    u_N-1 \\
    u_N
    \end{bmatrix}}_{x} = 
    \underbrace{
    \begin{bmatrix}
    -1 \\
    -1 \\
    \vdots \\
    -1\\
    -1
    \end{bmatrix}}_b
\end{equation}
Here, $A$ is the coefficient matrix, $b$ is the right hand side vector and $x$ is the unknown solution vector, which can be exactly obtained by inverting the matrix $A$, i.e., $x = A^{-1}b$. But, we will obtain an approximation of $x$ using iterative algorithms. We will instantiate two classical iterative algorithms: Gauss--Seidel and Jacobi, with this example problem and apply Theorem~\ref{nec_suf} to prove convergence of the approximate solutions, obtained using these algorithms, to the exact solution.

\vspace{-1em}
\subsection{Gauss--Seidel method}\label{sec:gauss_seidel_model}
We next demonstrate the convergence of the Gauss--Seidel iteration on the example~(\ref{matrix_A}). We choose $P=1$. Thus, we have a symmetric tri-diagonal coefficient matrix of size $3 \times 3$, which we will denote as $A_{GS}$. 
To show that iterative system for the system $A_{GS} x = b$ converges, we need to show that $A_{GS}$ is positive definite by application of the theorem \lstinline{Gauss_Seidel_converges}, which we proved in Section~\ref{sec:gauss_seidel_generic}. In Coq, we prove that $A_{GS}$ is positive definite (by proving positivity of determinant of the principal minors of $A_{GS}$, as illustrated in Example~\ref{example_pd}) the following lemma statement
\begin{lstlisting}
Lemma Ah_pd: forall (h:R), (0<h) -> is_positive_definite (Ah 2%N h).
\end{lstlisting}
Proving that $A_{GS}$ is positive definite using the approach illustrated in Example~\ref{example_pd} for a generic $N$ is tedious and does not add much to our line of argument. Hence, we chose to do it for $A_{GS}$ of size $3 \times 3$. One can perform this exercise for any choice of $N$ by defining an algorithm for computing the determinant of principal minors 
of a matrix and  get the same result. The statement of convergence of Gauss--Seidel iteration method for the $3 \times 3$ matrix is stated in Coq as
\begin{lstlisting}
Theorem Gauss_seidel_Ah_converges: forall (b: 'cV[R]_3) (h:R), (0 < h) -> 
let A := (Ah 2%N h) in
let x:= (A^-1) *m b in
forall x0: 'cV[R]_3, is_lim_seq (fun k:nat => vec_norm ((X_m k.+1 x0 b (M_G A) (N_G A)) -  x)) 0.
\end{lstlisting}
This closes the proof of the convergence of the Gauss--Seidel iteration for the model problem.

\vspace{-0.5em}
\subsection{Jacobi method}\label{sec:jacobi_proof}
We next apply the Theorem~\ref{nec_suf} to prove convergence of the Jacobi iteration on the model problem~(\ref{matrix_A}). As discussed earlier, the iteration matrix for a Jacobi iteration method is $S_J = I - D^{-1}A_J$.
We choose $P=1$, thereby obtaining a $3 \times 3$ matrix system, like we did for the Gauss--Seidel iteration. 

\vspace{-1em }
\subsubsection*{Formalization in Coq: }
We prove the following theorem in Coq for convergence of the Jacobi method for the model problem
\begin{lstlisting}
Theorem Jacobi_converges: forall (b: 'cV[R]_3) (h:R), (0 < h) -> 
let A := (Ah 2h) in 
let x := (A^-1) *m b in  
forall x0: 'cV[R]_3, is_lim_seq (fun k:nat =>vec_norm ((X_m k.+1 x0 b (M_J 2 h) (N_J 2 h)) - x)) 0.
\end{lstlisting}
To prove \lstinline{jacobi_converges} using the Theorem~\ref{nec_suf}, we need to prove that the modulus of each of the eigenvalues of $S_J$ is less than $1$. We prove this using the following lemma statement in Coq
\begin{lstlisting}
Theorem eig_less_than_1: forall (n:nat) (i: 'I_n.+1) (h:R),
  (0 < h) -> (0 < n) -> (C_mod (lambda_J i n h) < 1). 
\end{lstlisting}
Since the iteration matrix $S_J$ is tri-diagonal, we can define a closed form expression for \lstinline{lambda_J} using the formula
\begin{equation}
    |\lambda_i(S_J)|= \left| 1+\frac{h^2}{2}\lambda_i(A_J) \right|= \left| \cos{\left(\frac{m\pi}{P+1}\right)}
    \right|
\end{equation}

\noindent A caveat in using \lstinline{lambda_J} explicitly is that the definition of eigenvalue, \lstinline{lambda} in Theorem~\ref{nec_suf} is based on the roots of the polynomials in the diagonal of the Smith Normal form of a matrix. However, we use the closed form expression for eigenvalue of the tridiagonal iteration matrix $S_J$ for our model problem. Ideally, we would like to prove that \lstinline{lambda_J} equals \lstinline{lambda}. Since the proof of this relation is a tangent to the line of argument we are making in this work, we assume that this relation holds in this work, which we state formally in Coq as the following hypothesis
\begin{lstlisting}
Hypothesis Lambda_eq: forall (n:nat) (h:R) (i: 'I_n.+1), 
  let S_mat := RtoC_mat (- ( (M_J n h)^-1 *m (N_J n h) )) in
  lambda S_mat i = lambda_J i n h.
\end{lstlisting}
This closes the proof of iterative convergence for Jacobi iteration on the model problem.         

\vspace{-1em}
\section{Related work}
\vspace{-0.5em}
A number of works have recently emerged in the area of formalization of numerical analysis. This has been facilitated by advancements in automatic and interactive theorem proving~\cite{boldo2015coquelicot,garillot2009packaging}. Some notable works in the formalization of numerical analysis include the formalization of Kantorovich theorem for convergence of Newton methods~\cite{pasca2010formal}, the formalization of the matrix canonical forms~\cite{cano:hal-00779376}, and the formalization of the Perron-Frobenius theorem in Isabelle/HOL~\cite{thiemann2021perron} for determining the growth rate of $A^n$ for small matrices $A$. 
Boldo~et~al~ \cite{boldo2010formal,boldo2014trusting,boldo2013wave} proved consistency, stability and convergence of a second-order centered scheme for the wave equation. Tekriwal et al.~\cite{10.1007/978-3-030-76384-8_20} formalized the Lax equivalence theorem to guarantee convergence of a generic class of finite difference schemes. Besides Coq, numerical analysis of ordinary differential equations (ODEs) has also been done in Isabelle/ HOL~\cite{immler2012numerical}. Immler et al.~\cite{immler2016flow,immler2019flow} formalized flows, Poincar\'e map of dynamical systems, and verified rigorous bounds on numerical algorithms in Isabelle/HOL. In~\cite{10.1007/978-3-319-06200-6_9}, Immler formalized a functional algorithm that computes enclosures of solutions of ODEs in Isabelle/HOL. Deniz~et~al. formally analyze~\cite{deniz2022formalization} the problem of heat conduction in Isabelle/HOL, which assumes a closed form solution in a $1$-D domain and is thus specific to the problem. 

Since most of the existing formalization on differential equations assume either a closed form solution or use integration schemes which are problem specific, these works do not provide a generalized framework for solving differential equations numerically. Our work however addresses the issue of generalizability by formalizing a framework for solving a linear system iteratively, which is the approach followed by most general purpose linear solvers like ODEPACK~\cite{hindmarsh1983odepack}, SUNDIALS~\cite{hindmarsh2005sundials}. 

\vspace{-1em}
\section{Conclusion and Future work}
\vspace{-1em}
In this work we formalized a \emph{generalized theorem} about convergence of the solutions of an iterative algorithm to the \emph{true numerical solution (direct solution)}. 
In this process, we clarify various details in the proof of convergence, which are missing in the classical numerical analysis literature. We then instantiate this \emph{generalized} theorem to two classical iterative methods -- the Gauss--Seidel method, and the Jacobi method on a model problem. Since it is cumbersome to compute the eigenvalues of a generic matrix system, and verify that its magnitude is less than $1$, we provide a much easier check for convergence, especially for the Gauss--Seidel method, called the \emph{Reich theorem}~\cite{reich1949convergence}, which relies on the structure of matrix $A$. By composing the proof of the Reich theorem~\cite{reich1949convergence} with the main iterative convergence theorem, we show convergence of the Gauss--Seidel method for this matrix $A$. Thus, our approach is \emph{modular}, and can be extended to prove convergence of the Gauss--Seidel method on any desired matrix structure for which one can formalize conditions similar to the Reich theorem. 
During our fomalization, we develop a library in Coq to deal with complex vectors and matrices. We defined absolute values of complex numbers, common properties of complex conjugates and operations on conjugate matrices and vectors. This development leverages the existing formalization~\cite{cohen:hal-00671809,cohen:inria-00593738} of complex numbers and matrices in \mathcomp.~ The overall length of the Coq code and proofs is about 8.5k lines of code. It took us about 8 person-months of full time work for the entire formalization.

This work could be extended to develop a generalized end-to-end framework for verification of stationary iterative methods, which includes floating-point error analysis for a concrete implementation of the algorithm in C. In a related work by Tekriwal~et~al.~\cite{cicmpaper}, the authors provide an end-to-end correctness proof for a concrete implementation of the Jacobi iteration, which is an instance of stationary iterative methods. They use a bound on the real iterative solution for the Jacobi method, which relies on a proof of convergence of this real iterative solution to the ``true'' numerical solution. Thus, a crucial step in extending this result to a generic stationary iterative method, would be to use a \emph{generalized} proof of convergence for a stationary iterative method. The results from this paper can thus be directly used in extending the results from~\cite{cicmpaper} to a concrete implementation of a generic stationary iterative method.

This work could also be extended to verify solutions of non-linear systems. Most physical systems behave non-linearly, and the analysis of these non-linear systems is usually done by linearlizing it around an optimal trajectory.
We also plan on extending this work to analyze more practical class of iterative methods called the \emph{Krylov subspace methods}~\cite{saad2003iterative}, which use stationary iterative methods like Jacobi methods as preconditioners. Thus, our work is foundational in the analysis of these complicated and practical iterative solvers.

\vspace{-1.5em}

\bibliographystyle{splncs04}
\bibliography{references}

\newpage
\appendix

\section{Proof of Lemma~\ref{lem:subproof_1}}\label{thm:lem1_proof}
\begin{proof}[Proof of Lemma~\ref{lem:subproof_1}]
(\textbf{Necessity}):
We need to prove that
\begin{align*}
    &\lim_{k \to \infty} ||(-M^{-1}N)^k|| = 0 \implies (\forall x_0,  \lim_{k \to \infty} ||(-M^{-1}N)^k (x_0 - x)|| = 0)
\end{align*}
Given $x_0$,
\begin{align}
    &\lim_{k \to \infty} ||(-M^{-1}N)^{k} (x_{o}-x)|| \nonumber \\
    &\leq \lim_{k \to \infty}||(-M^{-1}N)^k||\,||x_{o}-x||;\quad [~||Ax||\leq ||A||~||x||~]\nonumber \\
    &= \left(\lim_{k \to \infty}||(-M^{-1}N)^k|| \right) \left(\lim_{k \to \infty}||x_{o}-x|| \right) \nonumber \\
    &=0;  \quad [\text{ since, } \lim_{k \to \infty}||(-M^{-1}N)^k||=0~]\nonumber 
\end{align}
(\textbf{Sufficiency}):
We need to prove that
\begin{align*}
    &(\forall x_0,  \lim_{k \to \infty} ||(-M^{-1}N)^k (x_0 - x)|| = 0) \implies \lim_{k \to \infty}||(-M^{-1}N)^k||=0
\end{align*}
We start by unfolding the definition of the norm of the iterative matrix
\begin{equation}\label{eqn:matrix_norm}
    ||(-M^{-1}N)^k|| = \sup_{(x_0 - x) \neq 0} \frac{||(-M^{-1}N)^k (x_0 - x) ||}{||x_0 - x||}
\end{equation}
Therefore, we need to prove that 
\begin{equation}\label{eqn:limit_matrix_norm}
    \lim_{k \to \infty}
    \left(\sup_{(x_0 - x) \neq 0}
    \frac{||(-M^{-1}N)^k (x_0 - x) ||}{||x_0 - x||} \right) = 0
\end{equation}
We can prove~(\ref{eqn:limit_matrix_norm}) by choosing an upper bound for the matrix norm in~(\ref{eqn:matrix_norm}), proving that the limit of this upper bound converges to zero and then applying the sandwich theorem~\cite{book_real} for limits. We choose this upper bound as $\sum_{j < n} ||(-M^{-1}N)^k e_j||$, i.e.,
\begin{equation}\label{eqn:upper_bound}
  \sup_{(x_0 - x) \neq 0} \frac{||(-M^{-1}N)^k (x_0 - x) ||}{||x_0 - x||} \leq \sum_{j < n} ||(-M^{-1}N)^k e_j||
\end{equation}
where $e_j$ is the unit vector corresponding to a principal direction in the Cartesian coordinate system, i.e., 
$e_j = \textbf{1}_j$. The vector $\textbf{1}_j$ is a unit vector with the entry $1$ in the $j^{th}$ place and other entries in the vector being $0$.
Therefore, 
\begin{align*}
     &\lim_{k \to \infty}
    \left(\sup_{(x_0 - x) \neq 0}
    \frac{||(-M^{-1}N)^k (x_0 - x) ||}{||x_0 - x||} \right) 
    \leq 
    \lim_{k \to \infty}
    \sum_{j < n} ||(-M^{-1}N)^k e_j|| \\
    &\implies
     \lim_{k \to \infty}
    \left(\sup_{(x_0 - x) \neq 0}
    \frac{||(-M^{-1}N)^k (x_0 - x) ||}{||x_0 - x||} \right) \leq \sum_{j < n}
    \left(\lim_{k \to \infty}
    ||(-M^{-1}N)^k e_j||\right)
\end{align*}
We can quantify $x_0$ in the hypothesis with $x + e_j, \forall j, j < n$. Therefore, $\forall j, j< n$, we have from the hypothesis,
\begin{equation*}
    \lim_{k \to \infty}
    || (-M^{-1}N) e_j|| = 0
\end{equation*}
Thus, 
\begin{equation*}
    \sum_{j < n}
    \left(\lim_{k \to \infty}
    ||(-M^{-1}N)^k e_j||\right) = 0
\end{equation*}
We can then apply the sandwich theorem~\cite{book_real} for limits to prove that
\begin{equation}\label{eqn:limit_matrix_norm_0}
     \lim_{k \to \infty}
    \left(\sup_{(x_0 - x) \neq 0}
    \frac{||(-M^{-1}N)^k (x_0 - x) ||}{||x_0 - x||} \right) = 0
\end{equation}
\end{proof}
\noindent We justify the choice of the upper bound in~(\ref{eqn:upper_bound}) in the Appendix~\ref{app:upper_bound}. Note that this proof was not available in the literature, and this is an original contribution as a result of our formalization.

\section{Proof of the equation~(\ref{eqn:upper_bound})}\label{app:upper_bound}
\begin{proof}[Proof of equation~(\ref{eqn:upper_bound})]
We can decompose the vector $x_0 - x$ into its components along the principal axes in the Cartesian coordinate system as
\begin{equation*}
    x_0 - x = \sum_{j < n} (x_0 - x)_j e_j
\end{equation*}
Therefore, 
\begin{equation}\label{eqn:linear_compose}
   (-M^{-1}N)^k (x_0 - x) = \sum_{j < n} (-M^{-1}N)^k (x_0 - x)_j e_j
\end{equation}
By taking a vector norm on both sides of~(\ref{eqn:linear_compose}),
\begin{equation*}
    ||(-M^{-1}N)^k (x_0 - x)|| =
    ||\sum_{j < n} (-M^{-1}N)^k (x_0 - x)_j e_j||
\end{equation*}
Using the triangle inequality property of the vector norm, we get
\begin{equation}\label{eqn:ineq_norm}
   ||(-M^{-1}N)^k (x_0 - x)|| \leq 
   \sum_{j < n} ||(-M^{-1}N)^k (x_0 - x)_j e_j||
\end{equation}
Since $(x_0 - x) \neq 0 $, $||x_0 - x|| \neq 0$. Hence, dividing by $||x_0 - x||$ on both sides of~(\ref{eqn:ineq_norm}),
\begin{align*}
    \frac{||(-M^{-1}N)^k (x_0 - x)||}{||x_0 - x||} &\leq 
   \sum_{j < n} \frac{||(-M^{-1}N)^k (x_0 - x)_j e_j||} {||x_0 -x||}\\
   &\leq 
   \sum_{j < n}
   \frac{|(x_0 - x)_j|}{||x_0 -x||}
    ||(-M^{-1}N)^k e_j|| \\
    &\leq \sum_{j < n}
    ||(-M^{-1}N)^k e_j||
\end{align*}
\end{proof}
\noindent Here, we first use the absolute homogenity property ($||ax|| = |a| ||x||$, for any scalar $a$ and vector $x$) of the vector norm. Then we use the fact that
\begin{equation*}
    |(x_0- x)_j| \leq 
    ||x_0 - x|| 
\end{equation*}

\section{Proof of Lemma~\ref{lem:subproof_2}}\label{thm:lem2_proof}
\noindent The quantity $\rho(-M^{-1}N)$ is the spectral radius of the iteration matrix $S$, and is defined as
\begin{equation*}
    \rho(S) = \max_{i} \{~|~\lambda_i(S)~|~\}, \forall i = 0,\hdots, (n-1)
\end{equation*}
where $\lambda_i(S)$ is the $i^{th}$ eigenvalue of $S$. Therefore,
\begin{equation*}
    \rho(S) < 1 \iff (\forall i, ~i < n \implies |\lambda_i(S)|)
\end{equation*}
\begin{proof}[Proof of Lemma~\ref{lem:subproof_2}]
(\textbf{Sufficiency}): We need to prove that
\begin{equation*}
     \lim_{k \to \infty} ||(-M^{-1}N)^k||=0 \implies \rho(-M^{-1}N) < 1
\end{equation*}
Since $\rho(-M^{-1}N) =\max_{0 \leq i < n}|\lambda_i (-M^{-1}N)|$,
\begin{equation*}
  \lim_{k \to \infty} ||(-M^{-1}N)^k||=0 \implies 
    (\forall i, 0 \leq i < n, \; |\lambda_i|<1)
\end{equation*}
Applying,
\begin{equation*}
  \lim_{k \to \infty} |\lambda_{i}|^k =0 \implies |\lambda_i|<1, \; \forall i, 0 \leq i < n  
\end{equation*}
to the goal statement, we need to prove:
\begin{equation*}
    \lim_{k \to \infty} |\lambda_i|^k = 0, \;\forall i, \; 0 \leq i < n
\end{equation*} 
under the hypothesis $\lim_{k \to \infty} ||(-M^{-1}N)^k||=0$. We now use the definition of an eigensystem, i.e.,
\begin{equation*}
   (-M^{-1}N) v_{i} = \lambda_i v_i 
\end{equation*}
where $v_i$ is an eigenvector corresponding to the eigenvalue $\lambda_i$. Therefore,
\begin{align}
   \lim_{k \to \infty} |\lambda_i|^k = 0 &\implies 
    \lim_{k \to \infty} |\lambda_i|^k ||v_i|| = 0 \nonumber \\
    & \implies \lim_{k \to \infty} ||\lambda_i ^k v_i||=0 \nonumber\\ 
    &\quad [\text{since, }|\lambda_i|^k = |\lambda_i^k| \; \land\; |\lambda_i^k|~ ||v_i|| = ||\lambda_i ^k v_i||] 
    \nonumber \\
    & \implies \lim_{k \to \infty} ||(-M^{-1}N)^k v_{i}||=0 \nonumber \\
    &\quad [\text{since, }(-M^{-1}N)^k v_{i} = \lambda_i ^k v_i] \label{part2_suf_1}
\end{align}
But the compatibility relation for vector and matrix norms dictates,
\begin{equation}\label{part2_suf_2}
    0 \leq ||(-M^{-1}N)^k v_{i}|| \leq ||(-M^{-1}N)^k||~ ||v_i||
\end{equation}
Since, $||v_i|| \neq 0 $ by definition of an eigensystem,
\begin{equation}
    \lim_{k \to \infty} ||(-M^{-1}N)^k||= 0 \implies \lim_{k \to \infty} ||(-M^{-1}N)^k||~ ||v_i|| =0
\end{equation}
We can then apply the sandwich theorem~\cite{book_real} to prove that 
\begin{equation*}
    \lim_{k \to \infty} || (-M^{-1}N)^k v_i|| = 0
\end{equation*}
\noindent
(\textbf{Necessity}):
We need to prove that
\begin{equation*}
    \rho(-M^{-1}N) <1 \implies \lim_{k \to \infty} ||(-M^{-1}N)^k||=0.
\end{equation*}
To prove necessity, we obtain the Jordan decomposition of the iterative matrix, $S=(-M^{-1}N)$. From the Jordan normal form theorem~\cite{saad2003iterative}, we know that there exists $V, J \in \mathbb{C}^{n \times n}$, $V$ non-singular and $J$ block diagonal such that: 
\begin{equation}\label{jordan_form}
    S= V J V^{-1}
\end{equation} with 
\begin{equation*}
    J = 
    \begin{bmatrix}
    J_{m_1}(\lambda_1) & 0 & 0 & \hdots & 0\\
    0 & J_{m_2}(\lambda_2) & 0 & \hdots & 0 \\
    \vdots & \hdots & \ddots & \hdots & \vdots \\
    0 & \hdots & 0 & J_{m_{s-1}}(\lambda_{s-1}) & 0\\
    0 & \hdots & \hdots & 0 & J_{m_s}(\lambda_{s})
    \end{bmatrix}
\end{equation*}
\begin{equation*}
     J_{m_i}(\lambda_i) = 
    \begin{bmatrix}
    \lambda_i & 1 & 0 & \hdots & 0\\
    0 & \lambda_i & 1 & \hdots & 0 \\
    \vdots & \vdots & \ddots & \ddots & \vdots \\
    0 & 0 & \hdots & \lambda_i & 1\\
    0 & 0 & \hdots & 0 & \lambda_i
    \end{bmatrix}
\end{equation*}
where $J_{m_i}(\lambda_i)$ is the Jordan block corresponding to the eigenvalue $\lambda_i$.
Thus, $S^k = V J^k V^{-1}$. Since $J$ is block diagonal,
\begin{equation*}
    J^m = 
    \begin{bmatrix}
    J_{m_1}^k(\lambda_1) & 0 & 0 & \hdots & 0\\
    0 & J_{m_2}^k(\lambda_2) & 0 & \hdots & 0 \\
    \vdots & \hdots & \ddots & \hdots & \vdots \\
    0 & \hdots & 0 & J_{m_{s-1}}^k(\lambda_{s-1}) & 0\\
    0 & \hdots & \hdots & 0 & J_{m_s}^k(\lambda_{s})
    \end{bmatrix}
\end{equation*}
Now, a standard result on the $k^{th}$ power of a $m_i \times m_i$ Jordan block states that, for $k \geq m_i-1$
\begin{equation*}
    J_{m_i}^k(\lambda_i) =
    \begin{bmatrix}
    \lambda_i^k & {k \choose 1}\lambda_i^{k-1} & {k \choose 2}\lambda_i^{k-2} & \hdots & {k \choose {m_i-1}}\lambda_i ^{k - m_i +1} \\
    0 &  \lambda_i^k & {k \choose 1}\lambda_i^{k-1} & \hdots &
    {k \choose {m_i-2}}\lambda_i ^{k - m_i +2}\\
    \vdots & \vdots & \ddots & \ddots & \vdots\\
    0 & 0 & \hdots & \lambda_i^k & {k \choose 1}\lambda_i^{k-1} \\
    0 & 0 & \hdots & 0 & \lambda_i^k
    \end{bmatrix}
\end{equation*}
To prove that $\lim_{k \to \infty} ||S^k||=0$, we need to prove that $\lim_{k \to \infty} ||J^k|| = 0 $, since $V$ is non-singular, i.e. $||V|| \neq 0$. Since the 2-norm of a matrix is bounded above by the Frobenius matrix norm as:
$ ||J^k||_2 \leq ||J^k||_F
$,
we can prove $\lim_{k \to \infty} ||J^k||_2 =0 $ by proving  $\lim_{k \to \infty} ||J^k||_F = 0$. The Frobenius norm of matrix $J^k$ is defined as:
\begin{equation*}
    ||J^k||_F \overset{\Delta}{=} \sqrt{\sum_{i} \sum_{j} |(J^k)_{(i,j)}|^2}
\end{equation*}
Thus, we need to prove that
\begin{equation}\label{each_entry_zero}
    \forall i \; j, \lim_{k \to \infty} |(J^k)_{(i,j)}|^2 = 0
\end{equation}
The zero entries of the Jordan block diagonal matrix tend to zero trivially.
Since, $|\lambda_i| < 1$, the diagonal elements of the Jordan block, $J_{m_i}^k (\lambda_i)$ tend to zero. Proving that the off diagonal elements tend to zero is more difficult.

We have to prove that $\lim_{k \to \infty} {k \choose m} |\lambda_i |^{k-m} = 0, \; \forall m \leq m_i -1$. The function $f(k) = {k \choose m}$ increases as $k$ increases while the function $g(k) = |\lambda_i|^{k-m}$ decreases as $k$ increases since $0 < |\lambda_i| < 1$. Informally, it can be argued that $g(k)$ decreases at a faster rate than the function $f(k)$. Hence, the limit should tend to zero. But proving it formally in Coq was challenging. We need to first bound the function $f(k)$ with a function $h(k) = \frac{k^m} {m!}$ to obtain a sequence, $\frac{k^m}{m!}|\lambda_i|^{k-m}$, for which it would be easy to prove the limit. Therefore, we split the proof into proofs of two facts:
\begin{itemize}
    \item 
    \begin{equation}\label{choose_prove}
    {k \choose m} \leq \frac{k^m}{m!}
    \end{equation}
    \item
    \begin{equation}\label{sequence_prove}
        \lim_{k \to \infty} \frac{k^m}{m!}|\lambda_i|^{k-m}=0
    \end{equation}
\end{itemize}
\noindent We provide a proof of the inequality~(\ref{choose_prove}) in the Appendix~\ref{app:ineq_29}. 
In order to prove $\lim_{k \to \infty} \frac{k^m}{m!} |\lambda_i|^{k-m} = 0$, we use the ratio test for convergence of sequences. The formalization of ratio test has not yet been done in Coq to our knowledge. So, we formalized the ratio test since it provides an easier test for proving convergence of sequences as compared to first bounding the function with an easier function for which the convergence could be proved with the existing Coq libraries. The process of bounding the function $\Gamma (k) = f(k)g(k)$ with an easier function to prove convergence was challenging for us due to the behavior of $\Gamma (k)$. Using plotting tools like Wolfram plot or MATLAB plot, we observed that for 
$|\lambda_i| \leq 0.5$, the function $\Gamma (k)$ was monotonously decreasing, while for $0.5 < |\lambda_i| < 1 $, $\Gamma (k)$ increases first and then decreases. Moreover, the location of the maxima of $\Gamma (k) $ in the interval $0.5 < |\lambda_i| < 1$ depends of the number of iterations. Hence, bounding $\Gamma (k)$ with a monotonically decreasing function is challenging for $0.5 < |\lambda_i| < 1$. But we know that $\Gamma (k) $ decays eventually. For such scenarios, just comparing the terms $\Gamma (k+1) $ and $\Gamma (k)$ provides a much simpler test for the convergence of the sequence $\Gamma (k)$. This is where ratio test for the convergence of sequence comes to rescue. The ratio test for the convergence of sequences is stated as follows.

\begin{lemma}\label{ratio_test}
If $(a_n)$ is a sequence of positive real numbers such that $\lim_{n \to \infty} \frac{a_{n+1}}{a_n} = L$ and $L < 1$, then $(a_n)$ converges and $\lim_{n \to \infty} a_n = 0$.
\end{lemma}
\noindent We provide an informal proof of the Lemma~\ref{ratio_test} in the Appendix~\ref{app:ratio_test}.
In our case, the sequence $a_n = n^m |\lambda_i|^n$. Therefore, the ratio $\frac{a_{n+1}}{a_n} = \frac{(n+1)^m|\lambda_i|^{n+1}}{n^m |\lambda_i|^n}$. Therefore,
\begin{align*}
    \lim_{n \to \infty} \frac{a_{n+1}}{a_n} = \frac{(n+1)^m |\lambda_i|^{n+1}}{n^m |\lambda_i|^n}
    &= \lim_{n \to \infty}\left( 1+ \frac{1}{n}\right)^m |\lambda_i|  \\
   & = |\lambda_i| \left[\lim_{n \to \infty}\left( 1+ \frac{1}{n}\right)^m\right] \\
    &= |\lambda_i| \left[\lim_{n \to \infty}\left( 1+ \frac{1}{n}\right)\right]^m  \\
    &=|\lambda_i|
\end{align*}
Therefore, $L=|\lambda_i|$, which we know is less than $1$. Thus, we can apply Theorem~\ref{ratio_test} to prove that $\lim_{n \to \infty} n^m |\lambda_i|^n = 0$. Since $m!$ and $|\lambda_i|^m$ are non zero constants,
$ lim_{k \to \infty} \frac{k^m}{m!} |\lambda_i|^{k-m} = 0. $
\end{proof}

\section{Proof of the inequality~(\ref{choose_prove})}\label{app:ineq_29}
The inequality~(\ref{choose_prove}) follows the following proof
\begin{proof}[Proof of equation~(\ref{choose_prove})]
\begin{equation}
   {m \choose k} = \frac{m!}{(m-k)! k!} 
   = \frac{(m-k+1)(m-k+2)...(m-1)m}{k!}\label{choose_num}
\end{equation}
Since $(m-k+1) \leq m$,  $(m-k+2) \leq m$, and so on, 
the product of terms in the numerator of~(\ref{choose_num}) is bounded by $m^k$. Hence, 
${m \choose k} \leq 
    \frac{m^k}{k!}$
\end{proof}

\section{Proof of the Lemma~\ref{ratio_test}}\label{app:ratio_test}
\begin{proof}
 Let $(a_n)$ be a sequence of real numbers such that $\lim_{n \to \infty} \frac{a_{n+1}}{a_n} = L$. Since $(a_n)$ is a sequence of positive numbers, $0 \leq L$. By the Density of Real Numbers theorem, there exists a real $r$, such that $L < r < 1$.

Unfolding the definition of $\lim_{n \to \infty} \frac{a_{n+1}}{a_n} = L$, for $\epsilon = r- L > 0 $, $\exists N \in \mathbb{N}$, such that forall $n \geq N$, $\left| \frac{a_{n+1}}{a_n} - L \right| < \epsilon = r - L$,  which implies that,
\begin{align*}
    &-\epsilon < \frac{a_{n+1}}{a_n} - L < \epsilon \nonumber \\
    &\implies L- \epsilon < \frac{a_{n+1}}{a_n} < L + \epsilon \nonumber \\
    &\implies \frac{a_{n+1}}{a_n} < L + (r-L) = r
\end{align*}
Therefore, $\forall n, \; n \geq N$, $a_{n+1} < a_n r$.
This implies, $a_n r < a_{n-1}r^2$, and $a_{n-2}r^2 < a_{n-3}r^3, \hdots ,a_{N+1}r^{n-N} < a_N r^{n-N+1}$.
We thus obtain the following inequality:
\begin{equation*}
    a_{n+1} < a_n r < a_{n-1}r^2 < \hdots < a_{N+1}r^{n-N} < a_N r^{n-N+1}
\end{equation*}
Let $K= \frac{a_N}{r^N}$. Therefore $a_N r^{n-N+1} = K r^{n+1}$. Thus for $n \geq N$, we have that $a_n r < K r^{n+1}$, or rather $a_n < K r^n$. Since $0 < r < 1$, we have that $\lim_{n \to \infty} r^n = 0$. Then by the squeeze lemma, $\lim_{n \to \infty} a_n = 0.$ 
\end{proof}

\end{document}